\documentclass[12pt]{article}
\usepackage{mathrsfs}
\usepackage{amsmath, amssymb, amsthm, latexsym, amsfonts, epsfig, color,graphicx,subfigure,subcaption,mathtools,extarrows}
\usepackage{float,booktabs}
\usepackage{algorithm,algpseudocode}
\allowdisplaybreaks[4]

\usepackage{listings}
\usepackage[colorlinks,
linkcolor=red,
anchorcolor=blue,
citecolor=green
]{hyperref}

\begin{document}
\newcommand{\dd}{{\rm d}}
\newcommand{\s}{\sigma}
\renewcommand{\k}{\kappa}
\newcommand{\p}{\partial}
\newcommand{\D}{\Delta}
\newcommand{\om}{\omega}
\newcommand{\Om}{\Omega}
\newcommand{\Ga}{\Gamma}
\renewcommand{\phi}{\varphi}
\newcommand{\e}{\epsilon}
\renewcommand{\a}{\alpha}
\renewcommand{\b}{\beta}
\newcommand{\N}{{\mathbb N}}
\newcommand{\R}{{\mathbb R}}
\newcommand{\eps}{\varepsilon}
\newcommand{\EX}{{\mathbb{E}}}
\newcommand{\PX}{{\mathbb{P}}}

\newcommand{\cF}{{\cal F}}
\newcommand{\cG}{{\cal G}}
\newcommand{\cD}{{\cal D}}
\newcommand{\cO}{{\cal O}}
\newcommand{\cA}{{\cal A}}
\newcommand{\de}{\delta}

\newcommand{\grad}{\nabla}
\newcommand{\n}{\nabla}
\newcommand{\curl}{\nabla \times}
\newcommand{\dive}{\nabla \cdot}

\newcommand{\ddt}{\frac{d}{dt}}
\newcommand{\la}{{\lambda}}

\newtheorem{theorem}{Theorem}[section]
\newtheorem{definition}[theorem]{Definition}
\newtheorem{lemma}[theorem]{Lemma}
\newtheorem{corollary}[theorem]{Corollary}
\newtheorem{example}[theorem]{Example}
\newtheorem{proposition}[theorem]{Proposition}
\renewcommand{\abstractname}{\Large\textbf{Abstract}}
\newtheorem{remark}{Remark}[section]
\newtheorem{assumption}{Assumption}[section]
\newtheorem{notation}{Notation}[section]
\def\proof{\mbox {\it Proof.~}}
\makeatletter
\@addtoreset{equation}{section}
\makeatother
\renewcommand{\theequation}{\arabic{section}.\arabic{equation}}

\allowdisplaybreaks[4]
\makeatletter\def\theequation{\arabic{section}.\arabic{equation}}\makeatother

\title{The governing equation for probability density of transition path time}
\author{
	{  Yan Liu, Xu Sun\footnote{Corresponding author}}\hspace{2mm}
	\vspace{1mm}\vspace{1mm}\\
	{\it\small  School of Mathematics and Statistics},\\
	{\it\small Huazhong University of Science and Technology},
	{\it\small Wuhan 430074,  China}\\
	{\it\small E-mail: liuyan1021@hust.edu.cn,
		xsun@hust.edu.cn}\vspace{1mm}}
\date{\today }
\maketitle




\setcounter{secnumdepth}{5} \setcounter{tocdepth}{5}

\makeatletter
    \newcommand\figcaption{\def\@captype{figure}\caption}
    \newcommand\tabcaption{\def\@captype{table}\caption}
\makeatother

\begin{abstract}
Transition path time  is the time it takes a system to move from the initial state (State A) to a final state (State B) without recrossing back to State A in between. The transition path time  is probabilistic in nature, and considerable research has been devoted to studying its distribution. In this paper, we derive a governing integral equation that expresses the density of the transition path time in terms of the probability density of the first passage time, which is well understood. An efficient numerical method is provided to solve this governing equation. Several numerical examples are presented to validate the governing equation and to compare the densities of the transition path time and the first passage time.
\end{abstract}


\section{Introduction}\label{section1}
Over the past few years, considerable research has focused on the stochastic dynamics of state transitions in complex systems. The transition path time between two states, say A and B, is defined as the time the system takes to move directly from A to B without returning to A. This quantity is crucial for understanding the dynamic mechanism underlying state transitions, including protein conformational changes and the folding dynamics of DNA hairpins \cite{chung2012, chung2013,sharma2020}.

Consider the following stochastic dynamical system:
\begin{equation}\label{e1.1}
	\left\{\begin{aligned}
		&\dd X(t)=f(X(t))\dd t+g(X(t))\dd B(t),\\
		&X(0)=x_0,
	\end{aligned}\right.
\end{equation}
where $X(t)$ is an $\R$-valued stochastic process, $B(t)$ is an $\R$-valued Brownian motion defined on a probability space $(\Om,\cF,\PX)$, $f:\R\rightarrow\R,g:\R\rightarrow\R$, and $x_0\in\R$ is a constant.

Let $\tau_{x,x_0}$ denote the first passage time  of the process $X(t)$  to reach $x$ from the initial state $x_0$, i.e.,
 \begin{equation}\label{e1.2}
 	\tau_{x,x_0}=\inf\{t>0| X(t)=x\}.
 \end{equation}
It is well known that 
\begin{align}\label{e1.2a}
\PX(\tau_{x,x_0}< \infty)=1.
\end{align}
There are well established methods for computing the probability density of the first passage time for $X(t)$ given in \eqref{e1.1}, see \cite{redner2001,sun2024}, among others.
For the system in \eqref{e1.1}, let $T_{x,x_0}$ denote the transition path time from the initial state $x_0$ to the target state $x$, then
 \begin{align}\label{e1.3}
	T_{x,x_0} =\begin{cases} \inf\{\xi>0| X(\xi)=x; \forall s \in (0,\xi), X(s)\ne x_0\}, \quad\text{if} \ \exists \lambda \in (0,\infty), X(\lambda)=x,\\
	\infty,\ \text{otherwise}.
\end{cases}
\end{align}
Note that 
\begin{align}\label{e1.3a}
\mathbb{P}(T_{x,x_0}<\infty) < 1, 
\end{align}
which contrasts with \eqref{e1.2a}.
Inequality \eqref{e1.3a} is a consequence of the fact that the system can return to the initial state $x_0$ before reaching the target state $x$. In fact, $T_{x,x_0}$ is defined to have a finite value only for paths that do not revisit $x_0$ before arriving at $x$, as can be seen from \eqref{e1.3}. 

Throughout this paper, we make the following assumptions:
\begin{assumption}\label{H1}
	The functions $f$ and $g$ are $C^\infty(\R)$ functions whose derivatives of any order are bounded.
\end{assumption}
\begin{assumption}\label{H2}
	$\forall \e >0, |g(x)|\geq \e$, for all $x\in\R$.
\end{assumption}
Under the assumptions stated above, the probability density  of   the first passage time is well-defined, see \cite{pauwels1987} for details. In this paper, the governing equation for the transition path time is derived, expressing its probability density  in terms of  the first passage time density through an integral equation. Consequently, the existence of the transition path time density follows directly from that of the first passage time if there exists a solution to the integral equation. The distribution of the transition path time has been measured experimentally in various contexts (see, e.g., \cite{gladrow2019, neupane2016}). On the theoretical side, analytical expressions have been derived under specific conditions such as near-equilibrium dynamics and high energy barriers \cite{lavacchi2024, laleman2017}. In other approaches, the distribution is obtained by numerically solving the  Fokker–Planck equation with absorbing boundary conditions \cite{sharma2020, caraglio2020}. However, these methods are essentially based on approximating the transition path time by the first passage time under certain restrictive assumptions. As a result, general and direct tools for computing transition path time distributions remain undeveloped. To the best of the authors’ knowledge, no governing equation exists for the distribution of  the transition path time in stochastic dynamical systems given by \eqref{e1.1}. The main objective of this paper is to derive such a governing equation.
 
This paper is organized as follows. In Section \ref{section2}, we derive a governing equation that governs the  density of the transition path time. A numerical method for solving this  governing equation is developed in Section \ref{section3}. In Section \ref{section4}, several numerical examples are presented to illustrate the proposed method and verify its efficiency. Section \ref{section5} is the conclusion.

\section{Governing Equation}\label{section2}
In this section, we shall derive the governing equation for the probability density of the transition path time. 

 It follows from \eqref{e1.2} and \eqref{e1.3} that
\begin{equation}\label{e2.1}
	\{\tau_{x,x_0}=t\}=\{X(t)=x;  X(0)=x_0; \forall s\in (0,t), X(s)\neq x\}
\end{equation}
and
\begin{equation}\label{e2.2}
\begin{aligned}
    \{T_{x,x_0}=t\}=&\{X(t)=x; X(0)=x_0; \forall s\in(0,t), X(s)\ne x_0\  \text{and}\  X(s)\ne x\}\\
    =&\{X(t)=x; X(0)=x_0; \forall u\in(0,t), X(u)\ne x_0;\forall v\in (0,t),  X(v)\ne x\},
\end{aligned}
\end{equation}
respectively. Let $\theta_{x,x_0}(t)=p(\tau_{x,x_0}=t|X(0)=x_0)$ represent the  probability density of the first passage time  $\tau_{x,x_0}$ at time $t$  for the process  $X(t)$ defined in \eqref{e1.1} with the initial value $x_0$,  and $\gamma_{x,x_0}(t)=p (T_{x,x_0}=t|X(0)=x_0)$  represent the probability density of the transition path time  $T_{x,x_0}$.  Then it follows from  \eqref{e2.1} and \eqref{e2.2}, respectively, that
 \begin{equation}\label{e2.3}
 \begin{aligned}
     \theta_{x,x_0}(t)=&	p(\tau_{x,x_0}=t|X(0)=x_0)\\
     =&p(X(t)=x;\forall u\in (0,t), X(u) \ne x |X(0)=x_0)
 \end{aligned}
\end{equation}
and
\begin{equation}\label{e2.4}
\begin{aligned}
     \gamma_{x,x_0}(t)=&p (T_{x,x_0}=t|X(0)=x_0)\\
     =&p(X(t)=x; \forall s\in (0,t), X(s)\neq x_0\;\text{and}\;X(s)\ne x|X(0)=x_0).
\end{aligned}
\end{equation}

For convenience, let $p(x,t|x_0,s)$  denote  $p(X(t)=x|X(0)=x_0)$, the  transition probability density  of $X(t)$ under the condition $X(0)=x_0$. Note that  \eqref{e2.3} can be written as
\begin{equation}\label{ee2.5}
    \begin{aligned}
        &\theta_{x,x_0}(t)\\
        =&p(\tau_{x,x_0}=t|X(0)=x_0)\\
        =&\int_0^t p(\tau_{x,x_0}=t; s \text{\ is the last time to visit}\  x_0 |X(0)=x_0 )\dd s
    \end{aligned}
\end{equation}
In fact,
\begin{equation}\label{ee2.6}
    \left\{s \text{\ is the last time to visit}\  x_0\right\}=\left\{X(s)=x_0; \forall v\in(s,t), X(v)\neq x_0\right\}.
\end{equation}
 Substituting \eqref{e2.1} and \eqref{ee2.6} into \eqref{ee2.5}, we obtain
\begin{equation}\label{e2.5}
\begin{aligned}
    &\theta_{x,x_0}(t)\\
    =&\int _0^t p\Big(X(t)=x;\forall u\in(0,t), X(u)\neq x; X(s)=x_0; \forall v\in(s,t), X(v)\neq x_0|X(0)=x_0\Big)\, \dd s\\
    =& \int_0^t p\Big(X(t)=x;\forall u\in(s,t),X(u)\neq x;\forall v\in (s,t),X(v)\neq x_0; X(s)=x_0; \\
     &\qquad\quad \forall u\in(0,s), X(u)\neq x|X(0)=x_0\Big)\, \dd s \\
     =&\int_0^t p\Big(X(t)=x;\forall u\in(s,t), X(u)\neq x; \forall v\in(s,t), X(v)\neq x_0\\
     &\quad\quad\quad \Big| X(s)=x_0; \; \forall u \in (0,s), X(u) \neq x; \; X(0)=x_0\Big) \\
     & \qquad\times p \big( X(s)=x_0; \; \forall u \in (0,s), X(u) \neq x \mid X(0)=x_0 \big)\, \dd s.\\
\end{aligned}
\end{equation}

Due to the time-homogeneous Markov  property of $X(t)$, the first factor of the integrand in the above integral can be written as
\begin{equation}\label{e2.6}
	\begin{aligned}
    &p(X(t)=x; \forall u\in(s,t), X(u)\neq x;\forall v\in(s,t), X(v)\neq x_0|X(s)=x_0; \\
    &\quad\ \forall u\in(0,s), X(u)\neq x; X(0)=x_0) \\
    =& p(X(t)=x; \forall u\in(s,t), X(u)\neq x; \forall v\in(s,t), X(v)\neq x_0|X(s)=x_0) \\
  =& p(X(t-s)=x;\forall u\in(0,t-s), X(u)\neq x;\forall v\in(0,t-s),X(v)\neq x_0|X(0)=x_0) \\
    = &\gamma_{x,x_0}(t-s).
	\end{aligned}
\end{equation}
The last ``$=$" in \eqref{e2.6} follows from \eqref{e2.4}. Note that the second factor of the integrand becomes
\begin{equation}\label{e2.7}
	\begin{aligned}
	 &p(X(s)=x_0;\forall u\in(0,s),X(u)\neq x|X(0)=x_0)\\
	 =&p(X(s)=x_0|X(0)=x_0)-p(X(s)=x_0; \exists u\in(0,s), X(u)= x|X(0)=x_0)\\
     =&p(x_0,s|x_0,0)-p(X(s)=x_0; \exists u\in(0,s), X(u)= x|X(0)=x_0).
	\end{aligned}
\end{equation}
The decomposition in \eqref{e2.7} follows from partitioning the event of returning to $x_0$ at time $s$ into those paths that avoid $x$  and those that hit $x$ before time $s$.

The second term in the right hand side of \eqref{e2.7} can be written as
\begin{align}
		&\quad p(X(s)=x_0; \exists u\in(0,s), X(u)= x|X(0)=x_0)\notag\\
        &=\int _0^s p(X(s)=x_0; \text{ the first visit to }\ x\ \text{occurs at time } u |X(0)=x_0)\dd u  \notag\\
        &=\int _0^s p(X(s)=x_0;\tau_{x,x_0}=u|X(0)=x_0)\dd u\notag\\
		&=\int _0^s p(X(s)=x_0|\tau_{x,x_0}=u; X(0)=x_0)p(\tau_{x,x_0}=u|X(0)=x_0)\dd u\label{e2.8}\\
        &\xlongequal{\text{by } \eqref{e2.1}}\int _0^s p(X(s)=x_0|X(u)= x; \forall  r\in(0,u), X(r)\neq x; X(0)=x_0) p(\tau_{x,x_0}=u|X(0)=x_0)\dd u\notag\\
        &=\int _0^s p(X(s)=x_0|X(u)= x; \forall  r\in(0,u), X(r)\neq x; X(0)=x_0) \theta_{x,x_0}(u)\dd u.\notag
\end{align}
	
By using the time-homogeneous Markov property of $X(t)$, the first factor of the above integrand becomes 
\begin{equation}\label{e2.10}
	\begin{aligned}
		&p(X(s)=x_0|X(u)= x;\forall  r\in(0,u),X(r)\neq x; X(0)=x_0)\\
		=&p(X(s)=x_0|X(u)= x)\\
		=&p(x_0,s-u|x,0).
	\end{aligned}
\end{equation}
Substitute \eqref{e2.10}  into \eqref{e2.8}, we get
\begin{equation}\label{e2.11}
	\begin{aligned}
		&p(X(s)=x_0;\exists u\in(0,s), X(u)= x|X(0)=x_0)\\
        =&\int _0^sp(x_0,s-u|x,0)\theta_{x,x_0}(u)\dd u.
	\end{aligned}
\end{equation}
Substitute \eqref{e2.11} into \eqref{e2.7}, we get 
\begin{equation}\label{ee2.11}
    \begin{aligned}
        &p(X(s)=x_0;\forall u\in(0,s),X(u)\neq x|X(0)=x_0)\\
        =&p(x_0,s|x_0,0)-\int _0^sp(x_0,s-u|x,0)\theta_{x,x_0}(u)\dd u.
    \end{aligned}
\end{equation}
Substituting \eqref{e2.6} and \eqref{ee2.11} into \eqref{e2.5} yields
\begin{equation}\label{e2.12}
	\begin{aligned}
		&\theta_{x,x_0}(t)\\
		=&\int_0^t \gamma_{x,x_0}(t-s) \left[p(x_0,s|x_0,0)-\int _0^s p(x_0,s-u|x,0)\theta_{x,x_0}(u)\dd u\right]\dd s,
	\end{aligned}
\end{equation}
which is the desired governing equation for the density of the transition path time. 

It is shown in \cite{sun2024} that $\theta_{x,x_0}(t)$ can be obtained by solving the following integral equation:
\begin{equation}\label{ee2.13}
    p(x,t|x_0,0)=\int_0^tp(x,s|x,0)\theta_{x,x_0}(t-s)\dd s.
\end{equation}
Once  $\theta_{x,x_0}(t),p(x_0,t|x_0,0)$ and $p(x_0,t|x,0)$ are known, $\gamma_{x,x_0}(t)$ can be obtained by solving the integral equation \eqref{e2.12}.
Taking the Laplace transform  at both sides of  \eqref{e2.12} and \eqref{ee2.13}, respectively, we get
\begin{equation}\label{e2.13}
\begin{aligned}
    	\overline{\theta}_{x,x_0}(\nu)=&\overline{\gamma}_{x,x_0}(\nu) \overline{p}_{x_0,x_0}(\nu)-\overline{\gamma}_{x,x_0}(\nu)\overline{p}_{x_0,x}(\nu)\overline{\theta}_{x,x_0}(\nu)
\end{aligned}
\end{equation}
and 
\begin{equation}\label{ee2.15}
    \overline{p}_{x,x_0}(\nu)=\overline{p}_{x,x}(\nu)\overline{\theta}_{x,x_0}(\nu),
\end{equation}
where
\begin{equation}\label{e2.14}
	\begin{aligned}
		\overline{p}_{x_2,x_1 }(\nu)=&\int_0^\infty p(x_2,t|x_1,0)e^{-\nu t}\dd t,\text{\ for\ } \forall x_1,x_2\in\R,\\
        \overline{\theta}_{x,x_0}(\nu)=&\int_0^\infty \theta_{x,x_0}(t) e^{-\nu t}\dd t,\\
        \overline{\gamma}_{x,x_0}(\nu)=&\int_0^\infty \gamma_{x,x_0}(t) e^{-\nu t}\dd t.\\
	\end{aligned}
\end{equation}
It follows from \eqref{e2.13} and \eqref{ee2.15} that
\begin{equation}\label{e2.15}
\begin{aligned}
    	\overline{\gamma}_{x,x_0}(\nu)=&\frac{\overline{\theta}_{x,x_0}(\nu)}{\overline{p}_{x_0,x_0}(\nu)-\overline{p}_{x_0,x}(\nu)\overline{\theta}_{x,x_0}(\nu)}\\
    =&\frac{ \overline{p}_{x,x_0}(\nu)}{\overline{p}_{x_0,x_0}(\nu)\overline{p}_{x,x}(\nu)-\overline{p}_{x_0,x}(\nu) \overline{p}_{x,x_0}(\nu)}.
\end{aligned}
\end{equation}

By applying the inverse Laplace transform to both sides of  \eqref{e2.15}, the solution for  ${\gamma}_{x,x_0}(t)$ can be obtained formally. However, performing this inversion  is often challenging, both analytically and numerically. 

We shall present an efficient numerical method to solve \eqref{e2.12} for the density of the transition path time.

\section{Numerical Method}\label{section3}
In order to solve $\gamma_{x,x_0}(t)$ numerically from the governing equation \eqref{e2.12}, we take the following steps:
\begin{enumerate}
	\item[(i)] Compute the transition probability densities $p(x_0,t|x_0,0),p(x_0,t|x,0),p(x,t|x_0,0)$ and $p(x,t|x,0)$  that appear in \eqref{e2.12} and \eqref{ee2.13}, and then compute $\theta_{x,x_0}(t)$, the probability density  of the first passage time, by solving \eqref{ee2.13}.
	\item[(ii)] Solve  $\gamma_{x,x_0}(t)$ from \eqref{e2.12}. 
\end{enumerate}

\subsection{The Transition Probability Density  and the First Passage Time Density}\label{section3.1}
Note that there are well established methods to compute the quantities $p(x_0,t|x_0,0)$ and $p(x_0,t|x,0)$ in step (i). We shall do this  computation using the method given in \cite{sun2024}. 

To proceed,  we divide the time interval $[0,T]$  into $N$ equal steps with size $\D t=\frac{T}{N}$, so that $t_k=k\D t$  for $k=0,1,\cdots,N$. Then, \eqref{e1.1} is discretized by the Euler-Maruyama scheme as
\begin{equation}\label{e3.1}
		\left\{
	\begin{aligned}
		&\hat{X}({t_k})=\hat{X}({t_{k-1})}+f(\hat{X}({t_{k-1})})\D t+g(\hat{X}({t_{k-1})})\triangle B(t_{k-1}),\text{ for }  k=1,2,\cdots,N,\\
		&\hat{X}(t_0)=x_0,
	\end{aligned}\right.
\end{equation}
where $\triangle B(t_{k-1})=B({t_k})-B(t_{k-1})$.

Let $\hat{p}(x_{k},t_{k}|x_{k-1},t_{k-1})$ denote the transition probability density of $\hat{X}(t_k)$ at $x_k$ conditioned on $\hat{X}(t_{k-1})=x_{k-1}$. It follows from \eqref{e3.1} that 
\begin{equation}\label{e3.2}
	\hat{p}(x_{k},t_{k}|x_{k-1},t_{k-1})=\frac{1}{\sqrt{2\pi g^2(x_{k-1})\D t}}\exp\left(-\frac{(x_{k}-x_{{k-1}}-f({x_{{k-1}}})\D t)^2}{2g^2(x_{k-1})\D t }\right).
\end{equation}
Under  Assumptions \ref{H1} and \ref{H2}, the probability density of $\hat{X}(t_k)$ converges to  that of  $X(t
_k)$ with first-order accuracy, as shown by the following lemma \cite[Corollary 2.1]{bally1996}. 
\begin{lemma}\label{lemma1}
		Suppose Assumptions \ref{H1} and \ref{H2} hold, then the probability density for $\hat{X}(t)$ converges to the probability density of $X(t)$, i.e. for $\forall j>i$ and $i,j\in\{0,1,2,\cdots,N\}$,
		\begin{equation}\label{e3.3}
			|\hat{p}(x, t_j|y,t_i)-p(x, t_j|y,t_i)|\leq C(T,x,y)\D t,
		\end{equation}
		where  $C(T,x,y)$ depends  on the time $T$ and the state variables $x,y$ such that $C(T,\cdot,\cdot)\in L^{\infty}(\R\times\R,\R^+)\bigcap L^1(\R\times\R,\R^+)$, i.e.
		\begin{equation}\label{e3.4}
			\forall x,y\in\R,\quad C(T,x,y)<\infty \ and \ \ \int_\R\int_\R C(T,x,y)\dd x\dd y<\infty.
		\end{equation}
\end{lemma}

For $t=\D t$, the transition probability density $\hat{p}(x,\D t|x_0,0)$  has an explicit expression given by  \eqref{e3.2}. For $t=t_k$, we propose to compute   $\hat{p}(x_{k},t_{k}|x_0,t_0)$  as follows.
\begin{equation}\label{e3.5}
	\begin{aligned}
		\hat{p}(x_{k},t_{k}|x_0,t_0)=&\overbrace{\int_\R \cdots\int_\R}^{k-1}\hat{p}(x_{k},t_{k}|x_{k-1},t_{k-1})\cdots\hat{p}(x_{2},t_{2}|x_{1},t_{1})\hat{p}(x_{1},t_{1}|x_0,t_0)\dd x_1\dd x_2\cdots\dd x_{k-1}.
	\end{aligned}
\end{equation}
Equation \eqref{e3.5} can be expressed recursively as
\begin{equation}\label{e3.6}
	\hat{p}(x_{k},t_{k}|x_0,t_0)=\int_\R\hat{p}(x_{k},t_{k}|x_{k-1},t_{k-1})\hat{p}(x_{k-1},t_{k-1}|x_{0},t_{0})\dd x_{k-1}.
\end{equation}
The above equation can be approximated as 
\begin{equation}\label{e3.7}
	\hat{p}(x_{k},t_{k}|x_0,t_0)=\int_{-D}^D\hat{p}(x_{k},t_{k}|x_{k-1},t_{k-1})\hat{p}(x_{k-1},t_{k-1}|x_{0},t_{0})\dd x_{k-1},
\end{equation}
where $D$ is a sufficiently large constant. We discretize the interval $[-D,D]$ into  $M+1$ nodes such that $-D=z_0< z_1<\cdots< z_{M}=D$ where $z_k=-D+k\D x,k=0,1,\cdots,M$ with $\D x=\frac{2D}{M}$. Equation \eqref{e3.7} can be expressed as
\begin{equation}\label{e3.9}
	\hat{p}(x_{k},t_{k}|x_0,t_0)=\sum_{i_{k-1}=0}^{M-1}\hat{p}(x_{k},t_{k}|z_{i_{k-1}},t_{k-1})\hat{p}(z_{i_{k-1}},t_{k-1}|x_{0},t_{0})\D x+\mathcal{O}(\D x),
\end{equation}
for $k=2,3,\cdots N$. It follows from \eqref{e3.2} and \eqref{e3.9}, for $k=1$,
\begin{equation}\label{e3.10}
	\hat{p}(x,t_{1}|x_0,t_0)=\frac{1}{\sqrt{2\pi g^2(x_0)\D t}}\exp\left(-\frac{(x-x_{0}-f({x_0})\D t )^2}{2g^2(x_{0})\D t}\right)
\end{equation}
and for $k=2,3,\cdots,N$,
\begin{equation}\label{e3.11}
	\begin{aligned}
		&\hat{p}(x_k,t_{k}|x_0,t_0)\\
        =&\D x\sum_{i_{k-1}=0}^{M-1}\frac{1}{\sqrt{2\pi g^2(z_{i_{k-1}})\D t}}\exp\left(-\frac{(x_k-z_{i_{k-1}}-f({z_{i_{k-1}}})\D t)^2}{2g^2(z_{i_{k-1}})\D t}\right)\hat{p}(z_{i_{k-1}},t_{k-1}|x_{0},t_{0})\\
        &+\mathcal{O}(\D x).
	\end{aligned}
\end{equation}

For simplicity,  define $\widetilde{p}$ as
\begin{equation}\label{ee3.11}
    \widetilde{p}(x,t_1|x_0,t_0)=\hat{p}(x,t_{1}|x_0,t_0)=\frac{1}{\sqrt{2\pi g^2(x_0)\D t}}\exp\left(-\frac{(x-x_{0}-f({x_0})\D t )^2}{2g^2(x_{0})\D t}\right)
\end{equation}
and 
\begin{equation}\label{ee3.12}
    \begin{aligned}
		&\widetilde{p}(x_k,t_{k}|x_0,t_0)\\
        =&\D x\sum_{i_{k-1}=0}^{M-1}\frac{1}{\sqrt{2\pi g^2(z_{i_{k-1}})\D t}}\exp\left(-\frac{(x_k-z_{i_{k-1}}-f({z_{i_{k-1}}})\D t)^2}{2g^2(z_{i_{k-1}})\D t}\right)\widetilde{p}(z_{i_{k-1}},t_{k-1}|x_{0},t_{0}),
	\end{aligned}
\end{equation}
for $k=2,3,\cdots,N$. Comparing \eqref{ee3.11} and \eqref{ee3.12} with \eqref{e3.9} and \eqref{e3.10}, we obtain
\begin{equation}
    \hat{p}(x,t_k|x_0,t_0)=\widetilde{p}(x,t_k|x_0,t_0)+\mathcal{O}(\D x), 
\end{equation}
which indicates $\widetilde{p}$  is a first-order approximation to  $\hat{p}$ with respect to $\D x$.
By Lemma \ref{lemma1},  
\begin{equation}\label{ee3.14}
    \begin{aligned}
        p(x,t_1|x_0,t_0)=\frac{1}{\sqrt{2\pi g^2(x_0)\D t}}\exp\left(-\frac{(x-x_{0}-f({x_0})\D t)^2}{2g^2(x_{0})\D t}\right)+\mathcal{O}(\D t).
    \end{aligned}
\end{equation}
and for $k=2,3,\cdots,N$,
\begin{equation}\label{e3.13}
	\begin{aligned}
		&{p}(x_{k},t_{k}|x_0,t_0)\\
		=&\sum_{i_{k-1}=0}^{M-1}\sum_{i_{k-2}=0}^{M-1}\cdots\sum_{i_1=0}^{M-1}\widetilde{p}(x_{k},t_{k}|z_{i_{k-1}},t_{k-1})\widetilde{p}(z_{i_{k-1}},t_{k-1}|z_{i_{k-2}},t_{k-2})\cdots \widetilde{p}(z_{i_{1}},t_{1}|x_{0},t_{0})(\D x)^{k-1}\\
		&+\mathcal{O}(\D x)+\mathcal{O}(\D t).
	\end{aligned}
\end{equation}
Consequently, $p$ is approximated  by  $\widetilde{p}$ with first-order accuracy with respect to  both $\Delta t$ and $\Delta x$. This establishes the numerical algorithm for computing the  density of the transition probability.

With the transition probability density numerically evaluated, we are now ready to compute  the probability density  $\theta_{x,x_0}(t)$ of  the first passage time by solving \eqref{ee2.13}. In fact, approximating each coefficient in \eqref{ee2.13} by \eqref{e3.13} and computing the integral in \eqref{ee2.13} by  the Euler scheme, we get \cite{sun2024}
\begin{equation}\label{ee3.16}
    \widetilde{p}(x,t_1|x_0,0)=\frac{\sqrt{\D t}}{\sqrt{2\pi}g(x_0)}\widetilde{\theta}_{x,x_0}(\D t)
\end{equation}
and for $k=2,3,\cdots,N$,
\begin{equation}\label{ee3.17}
    \begin{aligned}
        \widetilde{p}(x,t_k|x_0,0)=&\D t\sum_{i=2}^k\sqrt{i-1}(\sqrt{i}-\sqrt{i-2})\widetilde{p}(x,(i-1)\D t|x,0)\widetilde{\theta}_{x,x_0}((k-i+1)\D t)\\
        &+\frac{\sqrt{\D t}}{\sqrt{2\pi}g(x)}\widetilde{\theta}_{x,x_0}(k\D t).
    \end{aligned}
\end{equation}

According to \cite{sun2024}, the density of the first passage time, $\theta(t_k)$, in \eqref{ee2.13} is approximated  by $\widetilde{\theta}_{x,x_0}(t_k)$ for $k=1,2,\cdots, N$, which can be obtained by solving \eqref{ee3.16} and \eqref{ee3.17}.

\subsection{The Singularity of the  Governing Equation}
Before computing  the transition path time, it is essential to address  the singularity in the governing equation \eqref{e2.12}. Note that the right-hand side of \eqref{e2.12} involves a singular integral due to the fact that $\forall x_0\in\R$,
\begin{equation}\label{e3.14}
		\lim_{ t\rightarrow 0} p(x_0, t|x_0,0)
		\xlongequal{\text{by } \eqref{e3.3}}
        \lim_{t\rightarrow 0}\hat{p}(x_0,t|x_0,0)
		\xlongequal{\text{by } \eqref{e3.2}}\lim_{t\rightarrow 0}\frac{1}{\sqrt{2\pi  g^2(x_0) t}}\exp\left\{-\frac{f^2(x_0) t}{2 g^2(x_0) }\right\}
		=\infty.
\end{equation}
It is numerically challenging to solve \eqref{e2.12}  due to the singularity. To solve this problem,  let $s=\nu^2, \nu\geq 0$. Then  \eqref{e2.12} can be rewritten as 
\begin{equation}\label{e3.15}
	\begin{aligned}
		&\theta_{x,x_0}(t)\\
		=&2\int_0^{\sqrt{t}} \gamma_{x,x_0}(t-\nu^2) \left(p(x_0,\nu^2|x_0,0)-\int _0^{\nu^2} p(x_0,\nu^2-u|x,0)\theta_{x,x_0}(u)\dd u\right)\nu\dd \nu.
	\end{aligned}
\end{equation}
In this case,
\begin{equation}\label{e3.16}
	\begin{aligned}
		&\quad\lim_{\nu\rightarrow0}\nu\hat{p}(x_0,\nu^2|x_0,0)\\
         &\xlongequal{\text{by } \eqref{e3.2}}\lim_{\nu\rightarrow 0} \frac{\nu}{\sqrt{2\pi  g^2(x_0) \nu^2}}\exp\left\{-\frac{f^2(x_0) \nu^2}{2 g^2(x_0) }\right\}\\
		&=\frac{1}{\sqrt{2\pi  g^2(x_0) }}.
	\end{aligned}
\end{equation}
Note that 
\begin{equation}\label{e3.17}
	\begin{aligned}
		&\lim_{\nu\rightarrow0}\nu{p}(x_0,\nu^2|x_0,0)\\
		=&\lim_{\nu\rightarrow0}\nu\left({p}(x_0,\nu^2|x_0,0)-\hat{p}(x_0,\nu^2|x_0,0)\right)+\lim_{\nu\rightarrow0}\nu\hat{p}(x_0,\nu^2|x_0)\\
	\end{aligned}
\end{equation}
By Lemma \ref{lemma1},
\begin{equation}\label{e3.18}
	\lim_{\nu\rightarrow0}\nu\left({p}(x_0,\nu^2|x_0,0)-\hat{p}(x_0,\nu^2|x_0,0)\right)=0.
\end{equation}
Substituting \eqref{e3.16} and \eqref{e3.18} into \eqref{e3.17}, we get
\begin{equation}\label{e3.19}
	\lim_{\nu\rightarrow0}\nu{p}(x_0,\nu^2|x_0,0)=\frac{1}{\sqrt{2\pi  g^2(x_0) }}.
\end{equation}
Introduce
\begin{equation}\label{e3.20}
	\rho(x_0,\nu^2|x_0,0)=\left\{\begin{aligned}
		&\frac{1}{\sqrt{2\pi  g^2(x_0) }},\quad\  &\text{if}\ \nu=0,\\
		&\nu p(x_0,\nu^2|x_0,0),\ &\text{if}\ \nu>0.
	\end{aligned}\right.
\end{equation}
It follows from  \eqref{e3.19} and \eqref{e3.20} that $\rho(x_0,\nu^2|x_0,0)$ is continuous with respect to $\nu$. With the help of $\rho(x_0,\nu^2|x_0,0)$, the governing equation \eqref{e2.12} can be written as 
\begin{equation}\label{e3.21}
	\begin{aligned}
		&\theta_{x,x_0}(t)\\
		=&2\int_0^{\sqrt{t}} \gamma_{x,x_0}(t-\nu^2) \left(\rho(x_0,\nu^2|x_0,0)-\nu\int _0^{\nu^2} p(x_0,\nu^2-u|x,0)\theta_{x,x_0}(u)\dd u\right)\dd \nu.
	\end{aligned}
\end{equation}
Note that \eqref{e3.21} has no singularity and has the same solution as  \eqref{e2.12}.

\subsection{Solving the Governing  Equation}
Now we turn to solving \eqref{e3.21} numerically. 
Recall that the time interval $[0, t]$ is partitioned as $0=t_0<t_1<t_2<\cdots<t_{N-1}<t_N=t$ with  $t_k=k\D t,k=0,1,2,\cdots,N$.  Thus,  \eqref{e3.21} can be written as
\begin{equation}\label{e3.22}
	\begin{aligned}
		&\theta_{x,x_0}(t_k)\\
		=&2\sum_{i=1}^k\int_{\sqrt{t_{i-1}}}^{\sqrt{t_i}} \gamma_{x,x_0}(t_k-\nu^2) \left(\rho(x_0,\nu^2|x_0,0)-\nu\int _0^{\nu^2} {p}(x_0,\nu^2-u|x,0)\theta_{x,x_0}(u)\dd u\right)\dd \nu.
	\end{aligned}
\end{equation}
By the trapezoid  scheme, 
\begin{equation}\label{e3.23}
	\begin{aligned}
        &\int_{\sqrt{t_{i-1}}}^{\sqrt{t_i}} \gamma_{x,x_0}(t_k-\nu^2) \left(\rho(x_0,\nu^2|x_0,0)-\nu\int _0^{\nu^2} {p}(x_0,\nu^2-u|x,0)\theta_{x,x_0}(u)\dd u\right)\dd \nu\\
        =&\frac{1}{2}\left[\gamma_{x,x_0}(t_k-t_{i-1}) \left(\rho(x_0,t_{i-1}|x_0,0)-\sqrt{t_{i-1}}\int _0^{t_{i-1}} {p}(x_0,t_{i-1}-u|x,0)\theta_{x,x_0}(u)\dd u\right)\right.\\
    	&+\left.\gamma_{x,x_0}(t_k-t_i) \left(\rho(x_0,t_i|x_0,0)-\sqrt{t_i}\int _0^{t_i}  {p}(x_0,t_i-u|x,0)\theta_{x,x_0}(u)\dd u\right)\right](\sqrt{t_i}-\sqrt{t_{i-1}})\\
        &+E_i,
	\end{aligned}
\end{equation}
 where  $E_i$ is $\mathcal{O}(\sqrt{t_i}-\sqrt{t_{i-1}})^3$, i.e., there exists a constant  $C>0$ such that 
\begin{equation}\label{ee3.28}
    |E_i|\leq C(\sqrt{t_{i}}-\sqrt{t_{i-1}})^3=C(\sqrt{i}-\sqrt{i-1})^3(\D t)^{\frac{3}{2}}. 
\end{equation}
Substituting \eqref{e3.23} into \eqref{e3.22} yields
\begin{equation}\label{ee3.29}
	\begin{aligned}
		&\theta_{x,x_0}(t_k)\\
		=&\sum_{i=1}^{k}\left[\gamma_{x,x_0}(t_k-t_{i-1}) \left(\rho(x_0,t_{i-1}|x_0,0)-\sqrt{t_{i-1}}\int _0^{t_{i-1}} {p}(x_0,t_{i-1}-u|x,0)\theta_{x,x_0}(u)\dd u\right)\right.\\
    	&+\left.\gamma_{x,x_0}(t_k-t_i) \left(\rho(x_0,t_i|x_0,0)-\sqrt{t_i}\int _0^{t_i}  {p}(x_0,t_i-u|x,0)\theta_{x,x_0}(u)\dd u\right)\right](\sqrt{t_i}-\sqrt{t_{i-1}})\\
        &+2\sum_{i=1}^k E_i.
	\end{aligned}
\end{equation}
Using \eqref{ee3.28}, we get 
\begin{equation}\label{ee3.30}
\begin{aligned}
    &|\sum_{i=1}^k E_i|\\
    \leq &C \Big(1+(\sqrt{2}-1)^3+(\sqrt{3}-\sqrt{2})^3+\cdots+(\sqrt{k}-\sqrt{k-1})^3\Big)(\D t)^{\frac{3}{2}}\\    
    =&C\Big(1+\sum_{i=2}^k(\sqrt{i}-\sqrt{i-1})^3\Big)(\D t)^{\frac{3}{2}}\\    
    =&C\Big(1+\sum_{i=2}^k\big(\frac{1}{\sqrt{i}+\sqrt{i-1}}\big)^3\Big)(\D t)^{\frac{3}{2}}\\ 
    \leq &C\Big(1+\frac{1}{8}\sum_{i=1}^\infty i^{-\frac{3}{2}}\Big)(\D t)^{\frac{3}{2}}\\ 
    \leq & C'(\D t)^{\frac{3}{2}},
\end{aligned}
\end{equation}
where $C'=C\Big(1+\frac{1}{8}\sum_{i=1}^\infty i^{-\frac{3}{2}}\Big)>0$ is a constant independent of $k$ and $\D t$. This implies $\sum_{i=1}^k E_i$ is   $\mathcal{O}(\D t^{\frac{3}{2}})$.

Using the fact that $\sum_{i=1}^k E_i$ is   $\mathcal{O}(\D t^{\frac{3}{2}})$, as implied by  \eqref{ee3.30}, we obtain  for $k=1$, 
\begin{equation}\label{e3.27}
		\theta_{x,x_0}(t_1)=\rho(x_0,t_{0}|x_0,0)\gamma_{x,x_0}(t_1)\sqrt{\D t}+\mathcal{O}(\D t^\frac{3}{2}) 
\end{equation}
and for $k=2,3,\cdots, N,$ 
\begin{equation}\label{eq3.32}
	\begin{aligned}
		&\theta_{x,x_0}(t_k)\\
		=&\sum_{i=1}^{k}\left[\gamma_{x,x_0}(t_k-t_{i-1}) \left(\rho(x_0,t_{i-1}|x_0,0)-\sqrt{t_{i-1}}\int _0^{t_{i-1}} {p}(x_0,t_{i-1}-u|x,0)\theta_{x,x_0}(u)\dd u\right)\right.\\
    	&+\left.\gamma_{x,x_0}(t_k-t_i) \left(\rho(x_0,t_i|x_0,0)-\sqrt{t_i}\int _0^{t_i}  {p}(x_0,t_i-u|x,0)\theta_{x,x_0}(u)\dd u\right)\right](\sqrt{t_i}-\sqrt{t_{i-1}})\\
        &+\mathcal{O}(\D t^{\frac{3}{2}})\\
        =&\sum_{i=2}^{k}\gamma_{x,x_0}(t_k-t_{i-1})(\sqrt{t_i}-\sqrt{t_{i-1}}) \Big(\rho(x_0,t_{i-1}|x_0,0)\\
        &\qquad\qquad\qquad\qquad\qquad\qquad\qquad\quad\left.-\sqrt{t_{i-1}}\int _0^{t_{i-1}} {p}(x_0,t_{i-1}-u|x,0)\theta_{x,x_0}(u)\dd u\right)\\
        &+\sum_{i=1}^{k-1}\gamma_{x,x_0}(t_k-t_i) (\sqrt{t_i}-\sqrt{t_{i-1}})\left(\rho(x_0,t_i|x_0,0)-\sqrt{t_i}\int _0^{t_i}  {p}(x_0,t_i-u|x,0)\theta_{x,x_0}(u)\dd u\right)\\
        &+\gamma_{x,x_0}(t_k)\rho(x_0,t_0|x_0,0)\sqrt{\D t}+\mathcal{O}(\D t^{\frac{3}{2}})\\
        =&\sum_{i=1}^{k-1}\gamma_{x,x_0}(t_k-t_i) \left(\rho(x_0,t_i|x_0,0)-\sqrt{t_i}\int _0^{t_i}  {p}(x_0,t_i-u|x,0)\theta_{x,x_0}(u)\dd u\right)(\sqrt{t_{i+1}}-\sqrt{t_{i-1}})\\
        &+\gamma_{x,x_0}(t_k)\rho(x_0,t_0|x_0,0)\sqrt{\D t}+\mathcal{O}(\D t^{\frac{3}{2}})\\
	\end{aligned}
\end{equation}
To get \eqref{e3.27} and the second ``=" in \eqref{eq3.32}, we  use the fact that $\gamma_{x,x_0}(t_k-t_i)=0$ for $i=k$. 

By the trapezoid  scheme, the integral term in \eqref{eq3.32} becomes
\begin{equation}\label{eq3.33}
\begin{aligned}
      &\int _0^{t_i}  {p}(x_0,t_i-u|x,0)\theta_{x,x_0}(u)\dd u\\
      =&\sum_{j=1}^i \int_{t_{j-1}}^{t_j}  {p}(x_0,t_i-u|x,0)\theta_{x,x_0}(u)\dd u\\
      =& \frac{1}{2}\sum_{j=1}^i \Big[{p}(x_0,t_i-t_{j-1}|x,0)\theta_{x,x_0}(t_{j-1})+{p}(x_0,t_i-t_j|x,0)\theta_{x,x_0}(t_j)\Big]\D t+\mathcal{O}(\D t^2).
\end{aligned}
\end{equation}
Since $\theta_{x,x_0}(t_{j-1})=0$ for $j=1$ and $p(x_0,t_i-t_j|x,0)=0$ for $j=i$, \eqref{eq3.33} becomes 
\begin{equation}\label{eq3.34}
    \begin{aligned}
         &\int _0^{t_i}  {p}(x_0,t_i-u|x,0)\theta_{x,x_0}(u)\dd u\\
         =&\sum_{j=1}^{i-1} {p}(x_0,t_i-t_j|x,0)\theta_{x,x_0}(t_j)\D t+\mathcal{O}(\D t^2)
    \end{aligned}
\end{equation}
for $i=2,3,\cdots,k$. For case $i=1$, $\int _0^{t_1}  {p}(x_0,t_1-u|x,0)\theta_{x,x_0}(u)\dd u$ is of order $\D t^2$. Substituting \eqref{e3.20} and \eqref{eq3.34} into \eqref{e3.27} and \eqref{eq3.32}, respectively, we have for $k=1$,
\begin{equation}\label{eq3.35}
    	\theta_{x,x_0}(t_1)=\frac{\sqrt{\D t}}{\sqrt{2\pi g^2(x_0)}}\gamma_{x,x_0}(t_1)+\mathcal{O}(\D t^\frac{3}{2}),
\end{equation}
for $k=2$,
\begin{equation}\label{eq3.36}
    \theta_{x,x_0}(t_2)=\sqrt{{2}}\D t\gamma_{x,x_0}(t_1)p(x_0,\D t |x_0,0) +\frac{\sqrt{\D t}}{\sqrt{2\pi g^2(x_0)}}\gamma_{x,x_0}(t_2)+\mathcal{O}(\D t^{\frac{3}{2}})
\end{equation}
and for $k=3,4,\cdots,N$
\begin{equation}\label{eq3.37}
    \begin{aligned}
         &\theta_{x,x_0}(t_k)\\
        =&\D t\sum_{i=1}^{k-1}\sqrt{i} (\sqrt{{i+1}}-\sqrt{{i-1}})\gamma_{x,x_0}(t_k-t_i) \Big(p(x_0,t_i|x_0,0)\\
        &\qquad-\sqrt{t_i}\sum_{j=1}^{i-1} {p}(x_0,t_i-t_j|x,0)\theta_{x,x_0}(t_j)\D t\Big)+\frac{\sqrt{\D t}}{\sqrt{2\pi g^2(x_0)}}\gamma_{x,x_0}(t_k)+\mathcal{O}(\D t^{\frac{3}{2}}).
    \end{aligned}
\end{equation}

In \eqref{eq3.35}, \eqref{eq3.36}, and \eqref{eq3.37}, replacing $p$
with its approximation $\widetilde{p}$ (as justified by \eqref{ee3.14} and \eqref{e3.13}) and 
$\theta$ with its approximation $\widetilde{\theta}$ (as justified by \eqref{ee3.16} and \eqref{ee3.17}), and omitting the small error terms, allows these equations to be written in the following matrix form:
\begin{equation}\label{eq3.39}
	AX\approx b,
\end{equation}
where 
\begin{equation}
	X=\begin{bmatrix}
		\gamma_{x,x_0}(\D t)\\
		\gamma_{x,x_0}(2\D t)\\
		\vdots\\
		\gamma_{x,x_0}(N\D t)
	\end{bmatrix}, b=	\begin{bmatrix}
	\widetilde{\theta}_{x,x_0}(\D t)\\
	\widetilde{\theta}_{x,x_0}(2\D t)\\
	\vdots\\
	\widetilde{\theta}_{x,x_0}(N\D t)
	\end{bmatrix}
\end{equation}
and $A=(a_{ij})_{N\times N}$ with
\begin{equation}\label{eq3.40}
	a_{ij}=\left\{
	\begin{aligned}
		&0, &&{\rm for} \ i<j\\
		&\frac{\sqrt{\D t}}{\sqrt{2\pi g^2(x_0)}}, && {\rm for} \ i=j,\\
		&\sqrt{i-j}(\sqrt{i-j+1}-\sqrt{i-j-1})\D tI(({i-j})\D t), &&{\rm for} \ i>j.
	\end{aligned}\right.
\end{equation}
Here, 
\begin{equation}\label{eq3.38}
    I(t_i)=\left\{
    \begin{aligned}
    &\widetilde{p}(x_0,t_1|x_0,0),&&\text{for  } i=1,\\
     &\widetilde{p}(x_0,t_i|x_0,0)-\sqrt{t_i}\sum_{j=1}^{i-1} \widetilde{p}(x_0,t_i-t_j|x,0)\widetilde{\theta}_{x,x_0}(t_j)\D t &&\text{for } i=2,\cdots,N.       
    \end{aligned}\right.
\end{equation}

The transition path time  density vector $X$ is approximated by solving 
the following equation 
\begin{equation}
    A\widetilde{X}=b,
\end{equation}
where $\widetilde{X}=(\widetilde{\gamma}_{x,x_0}(t_1),\widetilde{\gamma}_{x,x_0}(t_2),\cdots,\widetilde{\gamma}_{x,x_0}(t_N))^T$ is  the approximate solution for $X$. The proposed algorithm for computing the probability density of the transition path time is summarized in Algorithm \ref{algorithm}.

Note that the probability density $\widetilde{\gamma}_{x,x_0}(t_k)$  contained in  $\widetilde{X}$  is unnormalized, as implied in \eqref{e1.3a}. The  normalized version  of the transition path time density, $\hat{\gamma}_{x,x_0}(t_k)$,  can be  obtained:
\begin{equation}\label{e3.42}
    \hat{\gamma}_{x,x_0}(t_k)=\frac{\widetilde{\gamma}_{x,x_0}(t_k)}{\sum_k \widetilde{\gamma}_{x,x_0}(t_k)}.
\end{equation}

\begin{algorithm}[htbp]
\caption{The proposed algorithm: computing the probability density of the transition path time for the system \eqref{e1.1} to reach the target state $x$ with the initial state $x_0$} 
\label{algorithm}

\begin{algorithmic}[1]

\Require $f,g,T,D,\D t$ and $\D x$. Here, $f,g$ are functions in \eqref{e1.1}, $T$ is the time, $D$ is the space domain, $\D t$ is the time step and $\D x$ is the space step

\Ensure $\hat{\gamma}_{x,x_0}(t_k),$ for $k=1,2,\cdots,N$, the normalized probability density of the transition path time for the system \eqref{e1.1}

\State set $N=\frac{T}{\D t}$ and $M=\frac{2D}{\D x}$

\For{$k = 0$ to $N$}
    \State set $t_k=k\D t$
\EndFor

\For{$k = 0$ to $M$}
    \State set $z_k=-D+k\D x$
\EndFor

\For{$k = 1$ to $N$}
    \State compute $\widetilde{p}(x_0,t_k|x_0,0), \widetilde{p}(x_0,t_k|x,0),\widetilde{p}(x,t_k|x_0,0)$ and $\widetilde{p}(x,t_k|x,0)$ by \eqref{ee3.11} and \eqref{ee3.12}
\EndFor

\For{$k = 1$ to $N$}
    \State  compute $\theta_{x,x_0}(t_k)$  by \eqref{ee3.16} and \eqref{ee3.17} with $\widetilde{p}(x,t_k|x_0,0)$ and $\widetilde{p}(x,t_k|x,0)$ obtained
\EndFor


\For{$i=1$ to $N$}
\For{$j=1$ to $N$}
\State compute $a_{ij}$ by \eqref{eq3.40}
\EndFor
\EndFor

\State get $A=(a_{ij})_{N\times N}$ with $a_{ij}$ obtained
\State get  $X$ by solving \eqref{eq3.39} 
\State get $\hat{\gamma}_{x,x_0}(t_k)$,  the normalized density by \eqref{e3.42} for $k=1,2,\cdots,N$
\State output $\hat{\gamma}_{x,x_0}(t_k)$ for $k=1,2,\cdots,N$
\end{algorithmic}
\end{algorithm}

\section{Examples}\label{section4}
In this section, three numerical examples are presented to illustrate the proposed algorithm.
\begin{example} \label{ex1}\rm 
Consider the following SDE:
	\begin{equation}\label{e4.1}
		\left\{\begin{aligned}
			&\dd X(t)=\dd B(t),\\
			&X(0)=x_0=0.
		\end{aligned}\right.
	\end{equation}
The solution to  \eqref{e4.1} corresponds to a standard Brownian motion.

The proposed method is applied to compute  the probability density  of the transition path time, where the time step is taken as $\Delta t = 1/1000$. Figure \ref{f11} compares the density by the proposed method with that by Monte Carlo simulations ($10^6$ runs). It can be observed from Figure \ref{f11} that the result obtained by the proposed method matches  well  with that by  Monte Carlo simulations.

Figure \ref{f12} compares the probability densities of the transition path time and the first passage time. The figure clearly shows that the transition path time has a more concentrated, left-skewed distribution compared to the first passage time. In fact, the transition path time is more probable to take much smaller values than the first passage time. 
 \begin{figure}[H]
    	\centering
    	\includegraphics[scale=0.8]{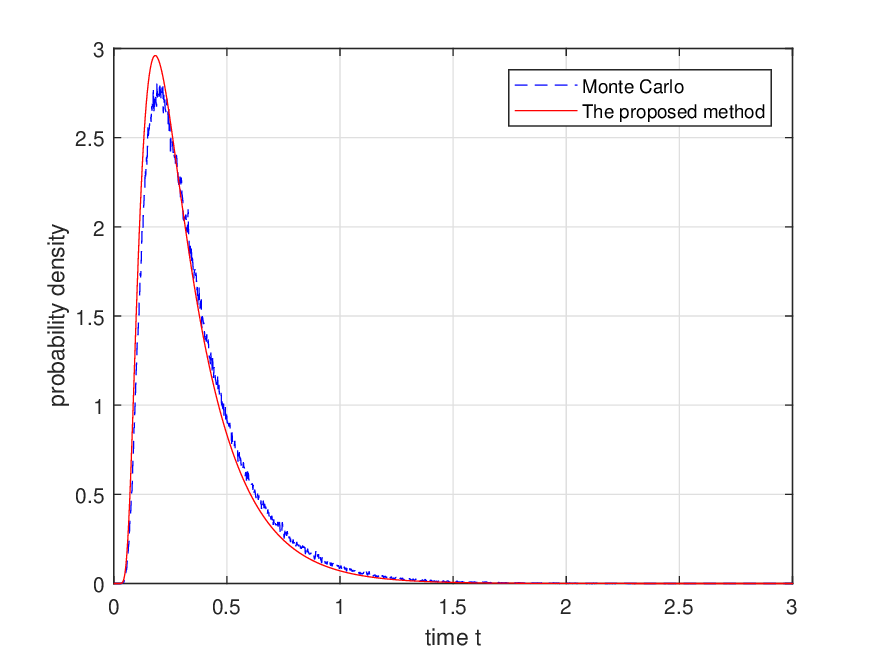}
    	\caption{Probability Density  for the System \eqref{e4.1} to Reach $x=1$ Starting from $x_0=0$: The Proposed Method vs Monte Carlo}
    	\label{f11}
    \end{figure}

     \begin{figure}[H]
    	\centering
    	\includegraphics[scale=0.8]{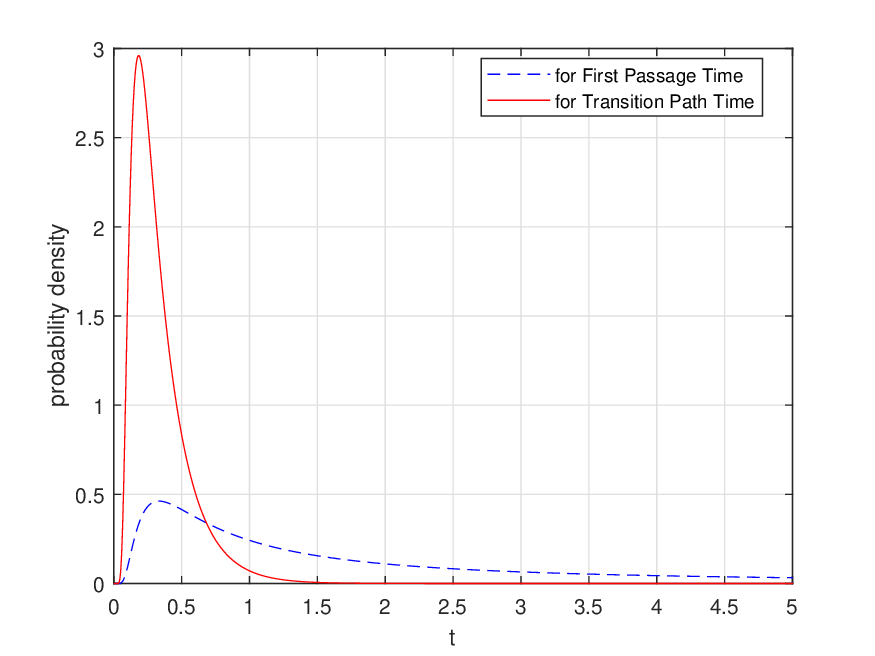}
    	\caption{Probability Density  for the System \eqref{e4.1} to Reach $x=1$ Starting from $x_0=0$: Transition Path Time vs First Passage Time}
    	\label{f12}
    \end{figure}


\end{example}
\begin{example}\label{ex2}\rm 
	Consider the following stochastic thermohaline circulation system with additive noise:
	\begin{equation}\label{e4.3}
		\left\{
		\begin{aligned}
			&\dd X(t)=(F-X(t)(1+\mu^2(1-X(t))^2))\dd t+\dd B(t),\\
			&X(0)=0.24,
		\end{aligned}\right.
	\end{equation}
    where $F=1.1$ represents a stochastic freshwater flux, $\mu^2=6.2$ is the diffusion ratio, and $B(t)$ represents  a standard Brownian motion. Without noise, the dynamical system \eqref{e4.3} possesses two stable equilibrium states  $x_1=0.24$ and $x_2=1.06$ \cite{sun2021}. Since the system can transition between these two states, we focus on the transition path time, which characterizes the time required for the state to move from $x_1=0.24$ to $x_2=1.06$.
	
	The probability density of the transition path time obtained by the proposed method is compared in Figure \ref{f21} with that from Monte Carlo simulations using $10^6$
  runs. The density of the transition path time is compared with that of  the first passage time in Figure \ref{f22}. In the computation, the time step is taken as $1/1000$ and the space step as $1/100$.  Similar to Example \ref{ex1}, Figures \ref{f21} and \ref{f22} show that the density obtained by the proposed method for the transition time agrees well with the   Monte Carlo simulation results. Moreover, the probability density of the transition path time is more likely to take smaller values than that of the first passage time.
  
    
    \begin{figure}[H]
    	\centering
    	\includegraphics[scale=0.8]{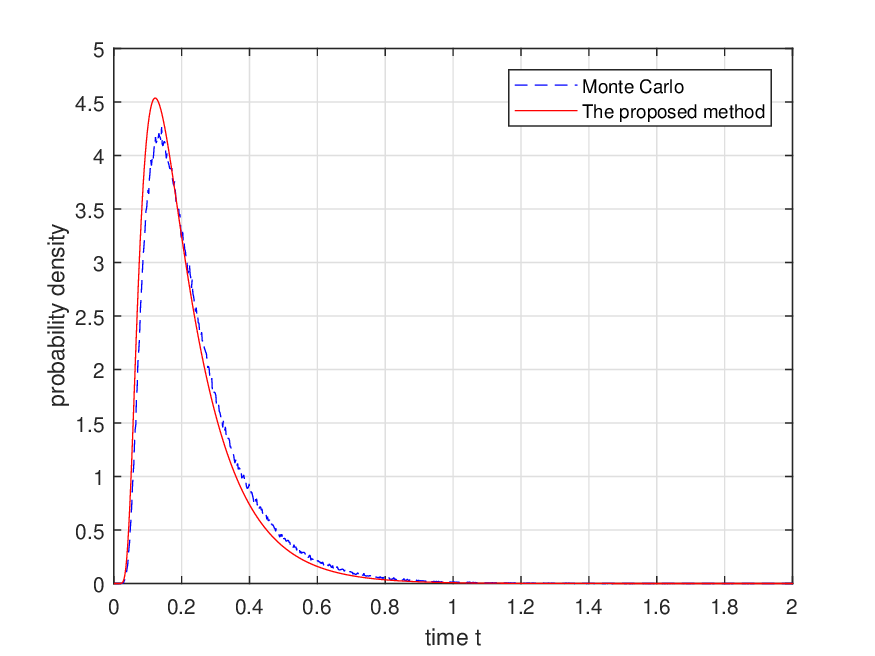}
    	\caption{Probability Density  for the System \eqref{e4.3} to Reach $x_2=1.06$ Starting from $x_1=0.24$: The Proposed Method vs Monte Carlo}
    	\label{f21}
    \end{figure}
    \begin{figure}[H]
    	\centering
    	\includegraphics[scale=0.8]{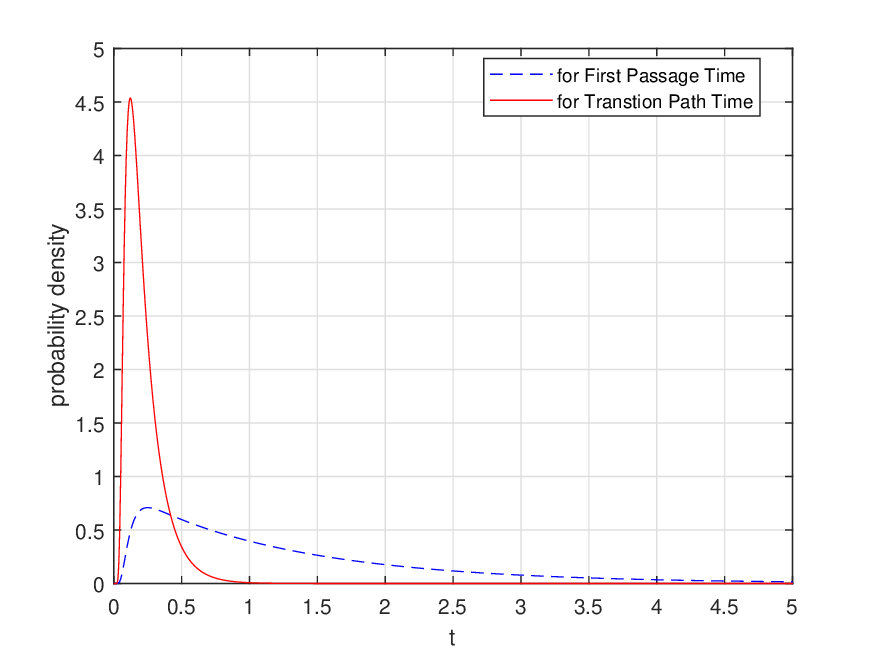}
    	\caption{Probability Density  for the System \eqref{e4.3} to Reach $x_2=1.06$ Starting from $x_1=0.24$: Transition Path Time vs First Passage Time}
    	\label{f22}
    \end{figure}

\end{example}

\begin{example}\label{ex3}\rm
		Consider the following stochastic thermohaline circulation system with multiplicative noise:
	\begin{equation}\label{e4.4}
		\left\{
		\begin{aligned}
			&\dd X(t)=(F-X(t)(1+\mu^2(1-X(t))^2))\dd t+\left(1+\frac{1}{2}\sin(X(t))\right)\dd B(t),\\
			&X(0)=0.24.
		\end{aligned}\right.
	\end{equation}
	This model differs form the one in Example \ref{ex2} by having multiplicative noise instead of additive noise. As in Example \ref{ex2}, we take $F=1.1$ and $\mu^2=6.2$. The probability density of the transition path time for $X(t)$ in system \eqref{e4.4} is numerically solved  with same parameters as in Example \ref{ex2}. Comparing  the densities by the proposed method and Monte Carlo simulation, as shown in Figure \ref{f31}, and densities for the transition path time and the first passage time, as shown in Figure \ref{f32}, we can get the similar observations as in Examples \ref{ex1} and \ref{ex2}.
    
By comparing Figure \ref{f31} here with Figure \ref{f21} from Examples \ref{ex2}, it is interesting to observe that the probability density curve becomes more concentrated around the peak and attains a higher peak value for cases with multiplicative noise than for those with additive noise.


     \begin{figure}[H]
    	\centering
    	\includegraphics[scale=0.8]{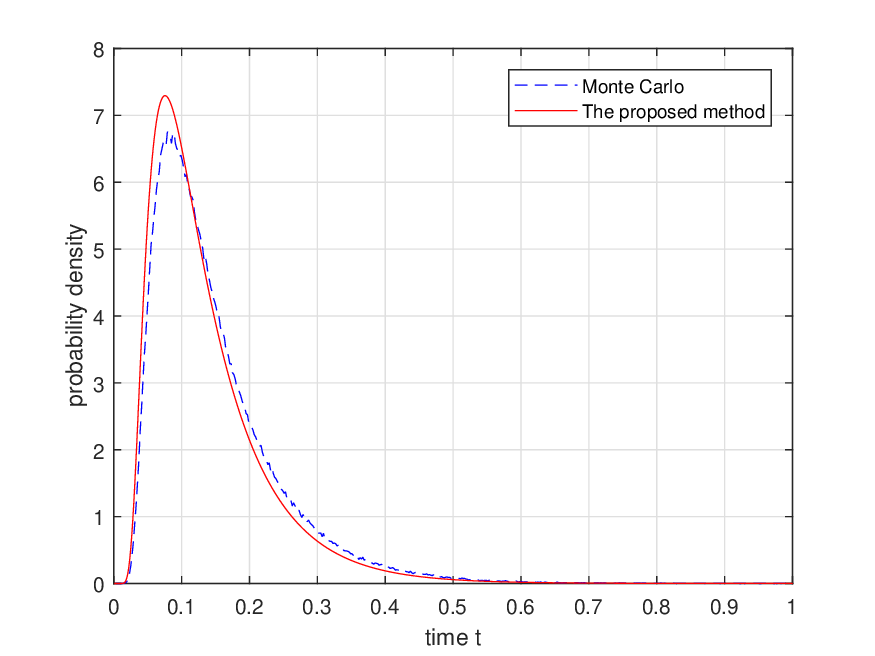}
    	\caption{Probability Density  for the System \eqref{e4.4} to Reach $x_2=1.06$ Starting from $x_1=0.24$: The Proposed Method vs Monte Carlo}
    	\label{f31}
    \end{figure}
         \begin{figure}[H]
    	\centering
    	\includegraphics[scale=0.8]{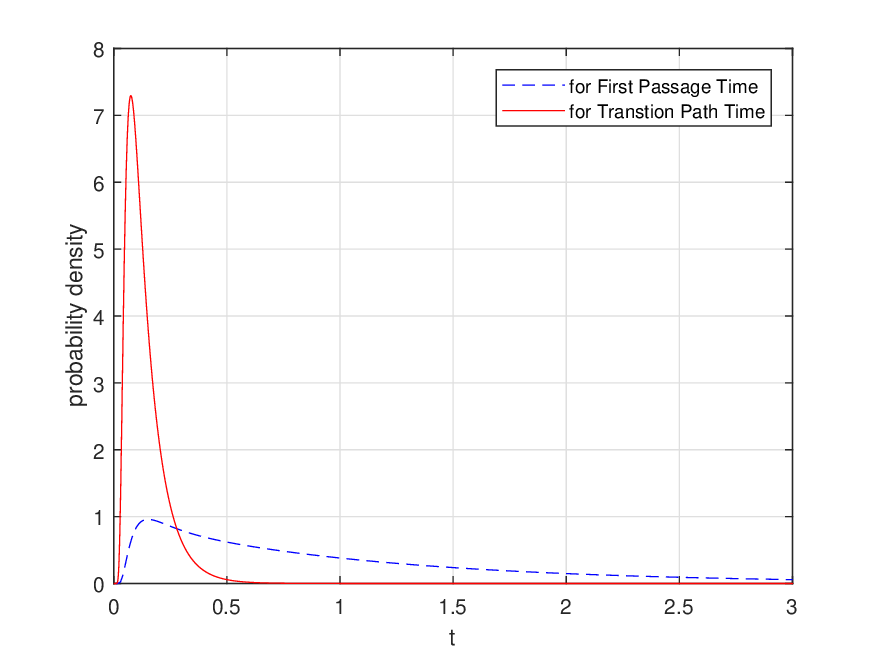}
    	\caption{Probability Density  for the System \eqref{e4.4} to Reach $x_2=1.06$ Starting from $x_1=0.24$: Transition Path Time vs First Passage Time}
    	\label{f32}
    \end{figure}
 \end{example}

\section{Conclusion}\label{section5}
In this paper, we derive a governing integral equation \eqref{e2.12}, for the probability density of the transition path time. An efficient numerical method, summarized in Algorithm \ref{algorithm}, is provided to solve this equation.

Note that the governing equation is derived only for the one-dimensional case. Generalization to higher dimensions is a nontrivial task and is left for future work.


\section*{Author Declaration}

\subsection*{Conflict of Interest }
 The authors have no conflicts to disclose.

 \section*{Data Availability } 
No datasets were generated or analysed during the
current study.





\newpage

\bibliography{ref}
\bibliographystyle{ieeetr}

\end{document}